# EL TEOREMA DE DIRICHLET, LA CONJETURA DE GOLDBACH Y ALGUNAS CONSECUENCIAS

**CAMPO ELIAS GONZALEZ PINEDA**

**cegp@utp.edu.co**

**RESUMEN.**

En este artículo, estudiamos y simplificamos el teorema de Dirichlet. Además hacemos una generalización del mismo, y deducimos de él la conjetura de Goldbach. También enunciamos y demostramos un teorema equivalente el postulado de Bertrand. Por último enunciamos unos resultados que se deducen el teorema de Goldbach-Elías .

# Abstract

In this article, we study and simplify Dirichlet's theorem. In addition we make a generalization of it, and deduce from it the Goldbach conjecture. We also state and prove an equivalent theorem Bertrand's postulate. Finally we enunciate some results that are deduced from the Goldbach-Elías theorem.

## CONTENIDO

1. **El teorema de Dirichlet:**

   Sean $a, b$ enteros primos relativos, la sucesión

$$p(n) = an + b$$

   Contiene infinitos números primos.

   **Ejemplo**: Existen infinitos números primos de la forma:

   - $pn + 1$ donde $p$ es un número primo.
   - $p^t n + 1$ con $p$ primo y $t$ un entero positivo.
   - $4pn + 1$con $p$ primo.
   - $(p_1 p_2 \dots p_t)n + 1$ donde cada $p_i$ es un número primo.
   - Si $p, q$ son números primos distintos $pn + q$.
   - $n + 1$. En este caso nótese que podemos hacer $n = am + b - 1$, por lo que $p(m) = am + b$.
   - En general existen infinitos números primos de la forma $2tn + 1$ e infinitos de la forma $2tn - 1$.

   **Nota1**.Es importante resaltar que si $p = an + b$ es un número primo entonces $a, b$ tienen que ser primos relativos. De la misma manera $n, b$ también tienen que ser primos relativos.

   **Nota2**. Aunque no es relevante la sucesión
   $p(n) = an + b$ es una sucesión aritmética. En efecto,

$$p(n) = a + b + a(n - 1), \qquad n \geq 1$$

   El siguiente resultado es de vital importancia en todo el proceso que seguiremos adelante.

   **Teorema2.**

   a) Sean $a, c$ impares positivos entonces,

   1. $\frac{a+c}{2}$ es par si y solo si $\frac{a-c}{2}$ es impar
   2. $\frac{a+c}{2}$ es impar si y solo si $\frac{a-c}{2}$ es par

   b) Si $a, c$ son enteros pares positivos $\frac{a+c}{2}$ y $\frac{a-c}{2}$ tienen la misma paridad.

**Demostración.**

a) Sean $a = 2t + 1, b = 2h + 1$ **entonces,**

$$\frac{a+c}{2} = \mathrm{t} + \mathrm{h} + 1, \text{ y } \frac{a-c}{2} = t - h$$

Veamos todas las posibilidades:

- $t\ y\ h$ pares entonces $t + h + 1$ es impar y $t - h$ par.
- $t\ y\ h$ impares entonces $t + h + 1$ es impar y $t - h$ par.
- $t$ par y $h$ impar entonces $t + h + 1$ es par y $t - h$ impar.

b) Sean $a = 2t, b = 2h$ entonces,

$$\frac{a+c}{2} = t + h, \qquad \frac{a-c}{2} = t - h$$

Si $t\ y\ h$ son pares entonces, $t + h$ es par y $t - h$ es par.

Si $t\ y\ h$ son impares entonces, $t + h$ es par y $t - h$ es par.

Si $t$ es par y $h$ impar entonces, $t + h$ es impar y $t - h$ es impar.

**2. Simplificando el teorema de Dirichlet.**

**Teorema3**. En el teorema de Dirichlet,

$$p(n) = an + b, \big((a, b)\big) = 1$$

Basta con que $a$ sea par, $b$ impar y $b < a$.

**Demostración.**

Es claro que $a, b$ no pueden ser simultáneamente pares, por tanto consideremos las siguientes posibilidades:

a) $a$ impar y b par. En este caso $n$ tiene que se impar, sea $n = 2t + 1$, luego,

$$p(t) = a(2t + 1) + b = 2at + a + b = a't + b'$$

Donde $a'$ es par y $b'$ es impar.

b) $a$ impar y $b$ impar. En este caso $n$ tiene que ser par, sea $n = 2t$, luego,

$$p(t) = a(2t) + b = 2at + b = a't + b$$

Donde $a'$ es par y $b$ es impar.

Así todas las posibilidades conducen a $a$ par y $b$ impar.

De otro lado si $a < b$

Por el algoritmo de la división,

$$b = at + r, \qquad 0 < r < a$$

Por lo que

$$p(n) = an + at + r = a(n + t) + r$$

Nótese que $r$ es impar.

**Ejemplo:**

1. La sucesión $p(t) = 2t + 1$ contiene todos los números primos impares.
2. Todo primo impar es de la forma $4n + 1$ o de la forma $4n - 1$ pero no de ambas a la vez.
3. Si $p(n) = an + b$ con $a$ par, $b < a$ impar, existe un impar $I$ tal que $b + I = a$, es decir $b = a - I$. Luego,

$$P(n) = an + a - I = a(n + 1) - I$$

Podemos concluir lo siguiente: Sean $a$ un entero par positivo e $I$ un impar positivo primos relativos. Existen infinitos números primos de la forma $p = an + I$ e infinitos de la forma $q = am - I$.

**3. Una consecuencia importante.**

Sea $a$ un par positivo y $b, c$ impares positivos primos relativos con $a$ y menores que $a$. (si $a = 2,\ b = c = 1$). Hagamos

$$p(n) = an + b, \qquad q(m) = am + c$$

Donde suponemos $p > q$. Promediando tenemos:

$$\frac{p(n) + q(m)}{2} = \frac{a}{2}(n + m) + \frac{b + c}{2} = A + B \quad (*)$$

Tenemos entonces las siguientes posibilidades.

1.

$$P = \frac{p(n)+q(m)}{2} \quad \text{par}$$

$$I = \frac{\boldsymbol{p(n)-q(m)}}{\mathbf{2}} \quad \text{impar}$$

De donde encontramos que

$$P + I = p(n)$$
$$P - I = q(m)$$

Y también

$$2P = p(n) + q(m)$$
$$2I = p(n) - q(m)$$

2.

$$I = \frac{p(n)+q(m)}{2} \quad \text{impar}$$
$$P = \frac{\boldsymbol{p(n)-q(m)}}{\boldsymbol{2}} \quad \text{par}$$

Luego,

$$I + P = p(n)$$

$$I - P = q(m)$$

De donde

$$2I = p(n) + q(m)$$

$$2P = p(n) - q(m)$$

Un examen más detallado puede hacerse $(*)$ analizando en lado derecho a $A = \frac{a}{2}(n+m), B = \frac{b+c}{2}$ pero esto se deja al lector.

Podemos entonces enunciar el siguiente resultado.

**Teorema4.**

1. Existen infinitos números pares de la forma $2P$ con $P$ par que se pueden escribir como la suma de dos números primos e infinitos de la forma $2I$ con $I$ impar que se pueden escribir como la resta de dos números primos.
2. Existen infinitos números pares de la forma $2I$ con $I$ impar que se pueden escribir como la suma de dos números primos e infinitos de la forma $2P$ con $P$ par que se pueden escribir como la resta de dos números primos.

Aunque el resultado se deduce del teorema de Dirichlet, también puede deducirse por mera observación. En efecto, como hay infinitos números primos $p$ basta hacer $2p = p + p$. También sabemos que existen infinitos números primos de la forma $p = 4k + 1$ e infinitos de la forma $q = 4l - 1$. Así tenemos infinitos pares de la forma $2h = p + q$.

De otro lado si $P$ es un número par positivo, vemos que existe un impar $I < P$ de tal manera que

$$2P = P + I + P - I = p + q$$

Donde claramente, $p, q$ son impares. Similarmente, si $I$es un impar positivo, existe un par $P < I$ de tal manera que

$$2I = I + P + I - P = p + q$$

Y una vez más , $p, q$ son impares.

La pregunta que nos hacemos es la siguiente: En la cuestión anterior, es posible elegir el impar $I$, (el par $P$) de tal manera que $p, q$ sean números primos?

Nuestro objetivo de aquí en adelante es dar respuesta a esta pregunta.

Comencemos analizando los siguientes casos.

1. Sean $p = 4n + 1, q = 4m + 1$ números primos y $p > q$ Entonces,

$$\frac{p+q}{2} = 2(n+m) + 1 = I$$
$$\frac{p-q}{2} = 2(n-m) = P$$

Por lo que

$$I + P = p$$
$$I - P = q$$

Y así,

$$2I = p + q$$
$$2P = p - q$$

2. Sean $p = 4n + 1, q = 4m - 1$ números primos y $p > q$ Entonces,

$$\frac{p+q}{2} = 2(n+m) = P$$
$$\frac{p-q}{2} = 2(n-m) + 1 = I$$

Por lo que

$$P + I = p$$
$$P - I = q$$

Y así,

$$2P = p + q$$
$$2I = p - q$$

3. Sean $p = 4n - 1, q = 4m - 1$ números primos y $p > q$ Entonces,

$$\frac{p+q}{2} = 2(n+m) - 1 = I$$
$$\frac{p-q}{2} = 2(n-m) = P$$

Por lo que

$$I + P = p$$
$$I - P = q$$

Y así,

$$2I = p + q$$
$$2P = p - q$$

Tenemos entonces el siguiente resultado.

**Teorema5**

1. Existen infinitos pares de la forma $2P$ para los cuales existe un impar $I < P$ de tal manera que

$$2P = P + I + P - I = p + q$$

Donde $p, q$ son números primos. El par $2I = p - q$.

2. Existen infinitos pares de la forma $2I$ para los cuales existe un impar $P < I$ de tal manera que

$$2I = I + P + I - P = p + q$$

Donde $p, q$ son números primos. El par $2P = p - q$.

### 3. Generalizando el teorema de Dirichlet

Recordemos la versión simplificada del teorema de Dirichlet.

**Teorema6 Elías- Dirichlet** Sea $t$ un entero positivo fijo. La sucesión

$$p(n) = 2tn + I, I \in \{1,3, \dots ,2t - 1\}$$

Con $2t$ e $I$ primos relativos contiene infinitos números primos.

Por un razonamiento claro, en el teorema anterior podemos escribir $p(n) = 2tn - I'$ ya que existe $I'$ tal que $I + I' = 2t$, es decir, $I = 2t - I'$. Así, $p(n) = 2t(n+1) - I'$.

Nótese que en el teorema de Dirichlet,

$$p(n) = an + b, \qquad \big((a,b)\big) = \big((n,b)\big) = 1$$

De manera similar, en el teorema anterior, $\big((2t, I)\big) = \big((n, I)\big) = 1$.

Veamos pues una generalización del teorema de Dirichlet.

**Teorema7 Elías** Para cada $n \geq 2$ existen $I, I' \in \{1,3, \dots 2tn-1\}$ de tal manera que:

1. $p(n) = 2tn + I$ es un número primo.
2. $p(n) = 2tn - I'$ es un número primo.

**Demostración.**

Vamos a demostrar este teorema por dos vías distintas, una es usando inducción matemática y la otra es utilizando el Postulado de Bertrand.

1. Hagamos el caso particular $t = 2$. Es decir, $p(n) = 4n + I$.
   Entonces, si $n = 2$, $p(2) = 8 + I$, $I \in \{1,3,5,7\}$
   Vemos que $I = 3$ o $I = 5$ hacen cierta la proposición. Supongamos que el resultado es cierto par $n$.Veamos para $n+1$.
   $$P(n+1) = 4(n+1) + I_1 = 4n + 4 + I_1, \;\; 4 + I_1 < 4(n+1)$$

Haciendo $I = 4 + I_1$y aplicando la hipótesis de inducción, encontramos que $p(n+1)$ es primo para algún $I < 4n < 4(n+1)$ y el resultado es cierto para todo $n \geq 2$.

De otro lado haciendo

$$p(n) = 2tn + I, t \geq 1$$

Vemos que $p(2) = 4t + I$es primo para algún $I < 4t < 2tn$ por lo que acabamos de demostrar. Si suponemos que

$$p(n) = 2tn + I$$

Es primo para algún $I < 2tn$ y hacemos

$$p(n+1) = 2t(n+1) + I_1 = 2tn + 2t + I_1$$

Es primo para algún $I = 2t + I_1 < 2t < 2tn$.

Justifiquemos ahora el resultado utilizando el Postulado de Bertrand.

Recordemos que el Postulado de Bertrand afirma que entre $n$ y $2n$ existe un primo $p$. Es decir, $n < p < 2n$.

Haciendo $n = 2tn$, existe un primo $p$ tal que

$$2tn < p < 4tn$$

Por tanto existe un impar $I$ de tal manera que $p = 2tn + I$.

Nótese que $I < 2tn$, ya que si $I \geq 2tn$, se tendría que $2tn + I \geq 4tn$, es decir, $p \geq 4tn$.

2. Demostremos esta parte utilizando el Postulado de Bertrand. Entre $tn$ y $2tn$ hay un primo $p$. Es decir, $tn < p < 2tn$. Sea $I$ impar tal que $p + I = 2tn$. Es decir,
$$p = 2tn - I.$$

El teorema de Elías y el Postulado de Bertrand son equivalentes, para mostrar esta afirmación, es necesario notar lo siguiente.

Sabemos que si $P, I$ son primos relativos, con $P$ par e $I$ impar existen infinitos primos de la forma

$$p(n) = In + P, P < I$$

Si $p(n)$ es primo la $n$ tiene que ser impar, es decir, $n = 2m + 1$. Luego,

$$p(m) = I(2m + 1) + P = 2Im + I_1$$

Donde $I_1 = I + P < 2I$.

De manera similar si

$$q(n) = In - P$$

Donde $I$ es impar, $P < I$ par, vemos que

$$-P < I \text{ es decir, } I - P < 2I$$

Para que $q(n)$ sea impar, $n$ tiene que ser impar. Si $n = 2m + 1$ tenemos,

$$q(2m + 1) = I(2m + 1) - P = 2Im + I_1, I_1 = I - P > 0.$$

Así tenemos el siguiente corolario.

**Corolario1** Sea $I$ un impar dado. Para cada $n \geq 5$ impar, existe $P, P' < In$ números pares tales que:

$$p = In + P, q = In - P'$$

Son números primos.

**Ejemplo.** Sea $I = 1$, entonces, $n + P$ y $n - P'$ con $n \geq 5$ son números primos para algún $P, P' \in \{0,2, \dots, n-1\}$. En particular si $n = 2m + 1$,

$$p(m) = 2m + 1 + P = 2m + I, I = P + 1$$

$$q(m) = 2m + 1 - P = 2m - I', I' = 1 - P'$$

Y además

$$I, I' \in \{1,3, \dots, 2m - 1\}$$

**Ejemplo.** Notar que existen infinitos primos de la forma $2n + I$, e infinitos de la forma $2n - I'$ donde $I, I' < 2n$ son impares.

## 4. Postulado de Bertrand

El postulado de Bertrand dice que: Entre $n$ y $2n$ hay un número primo. Es decir, existe un primo $p$ tal que

$$n < p < 2n$$

**Demostración.**

Por el teorema de Elías, sabemos que existe un primo $p$ para cada $n \geq 2$ tal que

$$p = 2n + I, \qquad I < 2n$$

Donde hemos elegido $t = 1$. Es claro entonces que

$$2n < 2n + I < 4n$$

De otro lado por el corolario anterior, si $p$ es un primo de la forma

$$p = In + P, P < In, \qquad P\ par$$

Entonces,

$$In < In + P < 2In$$

Haciendo $I = 1$ tenemos,

$$n < n + P < 2n$$

Esto completa la prueba.

En última podemos escribir

**Teorema8 Bertrand-Elías**

Para cada $n \geq 3$, existe un primo $p$ talque $n < p < 2n$ si y solo si para todo $n \geq 3$ existen impares $I, I'$ de tal manera que $2n + I, 2n - I'$ son números primos.

Es importante resaltar que si $p = 2n + I$ es un número primo, entonces $2n$ e $I$ son primos relativos. Sin embargo notemos lo siguiente.

**Ejemplo.** Suponga que $p = 2n - I$ es un número primo y que $\big((2n, I)\big) = l \neq 1$ entonces $l = p$.

En efecto, con las hipótesis dadas $2n = lk,\ I = k'l$ por lo que

$$p = 2n - I = l(k - k')$$

Como $p$ es primo y $l \neq 1$ se tiene que $k - k' = 1$ y $l = p$.

**Ejemplo.** Suponga que $\big((2n, I)\big) = q \neq 1$ entonces $p = 2n + I$ no es un número primo, ya que $q$ divide $p$.

En el primer ejemplo tenemos un caso especial, como caso particular $2(3) - 3 = 3$. Sin embargo siempre tenemos el siguiente resultado.

**Teorema9** Para todo $n \geq 2$, existe un impar $I < 2n$ de tal manera que $2n - I$ es un número primo y $\big((2n, I)\big) = 1$.

**Demostración.**

Es por inducción, hagamos $p(n) = 2n - I$. Para $n = 2$ tenemos,

$$p(2) = 4 - I, \qquad I \in \{1,3\}$$

Vemos que $I = 1$ hace cierta la proposición. Supongamos que el resultado se cumple para $n$. Ahora para $n + 1$,

$$p(n+1) = 2(n+1) - I_1 = 2n - (I_1 - 2)$$

Por hipótesis de inducción, existe $I = I_1 - 2 < 2n$ para el cual $2n - I$ es un número primo y $\big((2n, I)\big) = 1$. Hagamos

$$I_1 = I + 2$$

Es decir,

$$p(n+1) = 2(n+1) - I_1 = 2(n+1) - (I + 2), \qquad I + 2 < 2(n+1)$$

Ahora, sea $d = \big((2n + 2, I + 2)\big)$. Entonces, $d = 1$ o $d \neq 1$. Es claro que no puede ser $d \neq 1$ y por tanto $d = 1$.

**5. Resultado principal.**

Sabemos que para $n \geq 2$, existen impares $I_1,\ I_2 < 2n$ tales que

$$p(n) = 2n + I_1, \qquad q(n) = 2n - I_2$$

Son números primos. Queremos saber si es posible que $I_1 = I_2$.

Es decir, se trata de buscar $1 \leq I \leq 2n - 3$ de tal manera que

$$p(n) = 2n + I, \qquad q(n) = 2n - I \quad (*)$$

Sean números primos. Nótese que en $p + q = 4n$.

Para tal efecto, notemos que

$$1 \leq I \leq 2n - 3 \leftrightarrow 2n + 1 \leq 2n + I \leq 4n - 3$$

$$\leftrightarrow 2n + 1 \leq p \leq 4n - 3$$

De la misma manera

$$1 \leq I \leq 2n - 3 \leftrightarrow 3 - 2n \leq -I \leq -1$$

$$\leftrightarrow 3 \leq 2n - I \leq 2n - 1$$

$$\leftrightarrow 3 \leq q \leq 2n - 1$$

Sumando estas desigualdades obtenemos

$$2n + 4 \leq p + q \leq 6n - 4$$

Si suponemos que se cumple $(*)$ encontramos que

$$p + q = \frac{8n}{2} = 4n$$

Es decir, $p + q$ es el punto medio del intervalo $[2n + 4, 6n - 4]$. Notemos también que si se cumple $(*)$, $p = 2I + q$.

Dado el primo impar $3 \leq q \leq 2n - 1$ (notar que $q$ es de la forma $q = 2n - I$ ) veamos que existe $p$ primo de tal manera que $p = 2I + q$, con $I$ impar.

En efecto, como $2, q$ son primos relativos la sucesión

$$p = 2m + q$$

Contiene infinitos números primos(teorema de Dirichlet). Elijamos un impar $I$de tal manera que

$$p = 2I + q$$

Sea un número primo. Notemos que podemos escribir la ecuación anterior en la forma

$$p - I = I + q = 2l$$

Es decir,

$$p = 2l + I$$

$$q = 2l - I$$

Luego el $I$ que buscábamos existe al igual que el par $2l$ para nuestro objetivo deseado.

Es claro que

$$p + q = 4l = 4n \rightarrow l = n$$

Por un razonamiento similar tenemos lo siguiente:

Sea $q$ un número primo impar dado y consideremos la sucesión

$$p = 2m + q$$

Como $2, q$ son primos relativos esta sucesión contiene infinitos números primos. Elijamos un $H$ par de tal manera que

$$p = 2H + q$$

Sea un número primo. Es decir,

$$I = p - H = H + q,$$

$I$ impar, es decir

$$p = I + H$$

$$q = I - H$$

Es decir, $I$ tiene la propiedad de que existe $H$ par de tal manera que $p, q$ son números primos.

Recordemos que probamos que para todo impar $n \geq 5$ existen $H, H'$ tales que

$$p = in + H$$

$$q = in - H'$$

Donde $i$ es impar positivo y $H, H' \in \{0,2, \dots in - 1\}$. Si hacemos $i = 1$ tenemos

$$p = n + H$$

$$q = n - H'$$

Y así $H, H' \in \{0,2, \dots n - 1\}$. Busquemos $H$ tal que

a) $2 \leq H \leq n - 2 \leftrightarrow n + 2 \leq n + H \leq 2n - 2$
$\leftrightarrow n + 2 \leq p \leq 2n - 2$

b) $2 \leq H \leq n - 2 \leftrightarrow 2 - n \leq -H \leq -2$
$\leftrightarrow 2 \leq n - H \leq n - 2$
$\leftrightarrow 2 \leq q \leq n - 2$

Sumando encontramos

$$n + 4 \leq p + q \leq 3n - 4$$

Recordemos que

$$p = I + H$$

$$q = I - H$$

Sumando vemos $p + q = 2I$. Sumando las desigualdades obtenidas encontramos

$$p + q = \frac{n + 4 + 3n - 4}{2} = 2n = 2I \rightarrow n = I$$

Tenemos entonces los siguientes resultados.

**Teorema10 Elías-primos**

1. Para cada $n \geq 2$ exite un impar $I < 2n$ de tal manera que
$$p = 2n + I, q = 2n - I$$
Son números primos.
2. Para cada impar $n \geq 5$ existe un par $P < n$ de tal manera que
$$p = n + P, \qquad q = n - P$$
Son números primos.

**Corolario2. Goldbach-Elías** Todo entero para $P \geq 4$ se puede escribir como la suma de dos números primos.

**Nota importante:**

1. **Recordar que los primos $\boldsymbol{p}, \boldsymbol{q}$ están relacionados por la ecuación**
$$\boldsymbol{p = 2I + q}$$
2. En los intervalos $[2n + 4,\ 6n - 4], [n + 4, 3n - 4]$ El punto medio de cada intervalo es $p + q$. Si $n$ es par la suma de los primos es divisible por 4.
La longitud de cada intervalo es un número par:
$$l([2n + 4, 6n - 4]) = 4n - 8$$
$$l([n + 4, 3n - 4]) = 2n - 8$$
Nótese que los intervalos tienen la misma forma: $[a + 4, 3a - 4]$. La longitud es $2a - 8$ y el punto medio es $2a$.

3. Los intervalos los podemos escribir de la siguiente manera:
$$[2n + 4, 6n - 4] = [2n + 6, 6n + 2], \qquad n \geq 1$$
En este caso
$$p = 2n + 2 + I, q = 2n + 2 - I$$
También
$$[n + 4, 3n - 4] = [2n + 3,6n + 5], \qquad n \geq 1$$
En este caso
$$p = 2n - 1 + H, q = 2n - 1 - H$$

**6. Consecuencias**

En esta parte enunciamos unos resultados que se deducen del Teorema de Goldbach-Elías

1. Sea $P = 2n, n \geq 2$. Sabemos que existe un impar $I < 2n$ tal que

$$p = 2n + I \; q = 2n - I$$

Con $p, q$ números primos y que

$$\frac{p+q}{2} = 2n$$

Es fácil ver que los primos son de la forma

$$p = 4k + 1, q = 4l - 1$$

Por lo que

$$k + l = n$$

Así tenemos el corolario.

**Corolario Elias-2** Para cada $n \geq 2$ existen $k, l$ enteros positivos, tales que

$$k + l = n$$

Y

$$p = 4k + 1, q = 4l - 1$$

Son números primos.

Sea $I = 2n + 1,\ n \geq 1$. Sabemos que existe un par $P$ de tal manera que

$$p = I + P \; q = I - P$$

Con $p, q$ números primos y que

$$\frac{p+q}{2} = I = 2n + 1$$

Es fácil ver que los primos son de la forma

$$p = 4k + 1, q = 4l + 1$$

(o de la forma $p = 4k - 1, q = 4l - 1$)

Por lo que

$$k + l = n$$

El hecho de que los primos sean de la misma forma trae la consecuencia de cómo escribir el impar $I$. Que se puede escribir en la forma

$$I = 2n + 1 \;\; o \; I = 2n - 1.$$

Por ejemplo, $I = 19 = 2(9) + 1$. Notamos que

$$\begin{aligned} 9 &= 1 + 8 \\ &= 2 + 7 \\ &= 3 + 6 \\ &= 4 + 5 \end{aligned}$$

Y además:

$$4(1) + 1 = 5, \qquad 4(8) + 1 = 32$$

$$4(2) + 1 = 9, \qquad 4(7) + 1 = 29$$

$$4(3) + 1 = 13, \qquad 4(6) + 1 = 25$$

$$4(4) + 1 = 17, \qquad 4(5) + 1 = 21$$

Sin embargo, si escribimos $I = 2(10) - 1$ encontramos

$$\begin{aligned} 10 &= 1 + 9 \\ &= 2 + 8 \\ &= 3 + 7 \\ &= 4 + 6 \end{aligned}$$

Y así,

$$4(1) + 1 = 5, \qquad 4(9) + 1 = 37$$

$$4(2) + 1 = 9, \qquad 4(8) + 1 = 33$$

$$4(3) + 1 = 13, \qquad 4(7) + 1 = 29$$

$$4(4) + 1 = 17, \qquad 4(6) + 1 = 25$$

**Corolario-Elías3** Para todo impar $I = 2n + 1 = 2(n + 1) - 1, \; n \geq 2$ existen enteros positivos $k, l$

tales que

$$k + l = n \text{ o } k + l = n + 1$$

Y

$$p = 4k + 1, q = 4l + 1$$

Son números primos.

2. Sabemos que si $P$ es un número par mayor o igual a 4, existen primos $p, q$ tales que

$$P = p + q$$

Excluyendo el caso trivial $P = 4$, existe un impar $I$ de tal manera que

$$P + I = p$$
$$P - I = q$$

Donde $I < P$, $\left((P, I)\right) = 1$.

Restando las dos ecuaciones

$$2I = p - q$$

De la misma manera si $I$ es un impar mayor o igual a tres, existe un par $P$, tal que

$$I + P = p$$
$$I - P = q$$

con primos $p, q$ donde $P < I$ $\left((P, I)\right) = 1$.

Restando encontramos,

$$2P = p - q$$

Notar: $I \geq 5, I = 3$ es el caso trivial.

**Corolario4-Elías** Todo par $P \geq 4$ se puede escribir como la resta de dos números primos.

3. Por el teorema de Goldbach-Elías sabemos que todo par $2P \geq 4$ se puede escribir como la suma de dos números primos, digamos

$$2P = p + q$$

Si sumamos un primo impar $r$ tenemos,

$$2P + r = p + q + r$$

Por ejemplo si $r = 3$ tenemos la lista

$$6 = 3 + 3 \leftrightarrow 9 = 3 + 3 + 3$$

$$8 = 5 + 3 \leftrightarrow 11 = 5 + 3 + 3$$

Es decir, todo impar mayor o igual a 9 de la forma

$$2P + 3, \qquad P \geq 3$$

se puede escribir como la suma de tres primos impares.

$$2P + 3 = p + q + 3$$

Podría uno preguntarse, que impares mayores o iguales a 9 pueden escribirse como la suma de tres primos todos distintos.

**Corolario5** Todo impar mayor o igual a 9 se puede escribir como la suma de tres primos impares.

Notemos que si

$$2P + 3 = p + q + 3$$

Notemos varias cosas:

1. Si $p \geq 9, q < 9$ entonces,
$$2P + 3 = p_1 + p_2 + p_3 + q + 3$$
2. Ahora, si Si $p \geq 9, q \geq 9$ entonces
$$2P + 3 = p + q + 3 \geq 21$$

Y se tiene que

$$2P + 3 = p_1 + p_2 + p_3 + p_4 + p_5 + p_6 + 3$$

Una pregunta más general, sea $I \geq 3$ un impar, ¿cuál es el número máximo de números primos impares para el cual
$I = p_1 + p_2 + p_3 + \cdots + p_k$?

De estos $I$ cuáles tienen la propiedad de que $p_i \neq p_j$ para todo $i, j$?

Sin embargo notemos que si

$$I = I_1 + I_2 + \cdots + I_{2n+1}, n \geq 0$$

Donde cada $I_k \geq 9$ impar entonces $I$ se puede escribir como la suma de $3(2n+1)$ números primos. Notar que $I \geq 9(2n+1)$.

De otro lado si $I$es un impar y $p$ un número primo, entonces

$$I = pc + r, \quad 0 \leq r < c$$

Tenemos las siguientes posibilidades:

- Si $r < 9$ es primo, $I$ se puede escribir como la suma de $c + 1$ primos.
- Si $r$ es impar mayor o igual a 9, $I$ se puede escribir como la suma de $c + 3$ números primos.
- Si $r$ es un par mayor a 4, $I$ se puede escribir como la suma de $c + 2$ números primos.

Evidentemente hay más posibilidades pero estas se dejan al lector.

Notemos por ahora lo siguiente. Es claro que si $P \geq 4 \leftrightarrow 2P \geq 8$ con $P$ par, entonces

$$2P = P + P = p + q + p + q$$

De la misma manera si $I$ e un impar entonces.

$$2I = I + 1 + I - 1 = p + q + r + t$$

Donde $p, q, r, t$ son números primos y debe cumplirse que

$$\mathrm{I} - 1 \geq 4 \leftrightarrow \mathrm{I} \geq 5 \leftrightarrow 2\mathrm{I} \geq 10$$

En últimas tenemos,

**Teorema:** Todo par mayor o igual a 8 se puede escribir como la suma de cuatro primos impares.

**Observación.** Hagamos $I = 2a + 3, a \geq 3$. Podemos escribir

$$I = 2a + 3 + 2b - 2b \geq 9, \qquad b \geq 1, a \geq 3.$$

O los que es lo mismo

$$I = 2a + 2b + 3 \geq 2b + 9$$

Eligiendo $2b + 3$ un número primo, vemos que podemos reemplazar el primo 3 por cualquier otro primo impar. Podemos utilizar este resultado para encontrar números impares y decir cuantos primos impares es su suma.

Por ejemplo, si $b = 5$, el impar $2a + 13$ se puede escribir como la suma de cinco números primos impares.

1. **Primos gemelos**
   Diremos que dos primos son gemelos si son de la forma
   $$p = P + I, \qquad q = P - I$$

Es claro entonces que:

**Teorema** Existen infinitos primos gemelos.

**Referencias bibliográficas.**

[1].Una nota sobre la conjetura de Goldbach. Campo Elías González Pineda. Scientia Et Technica, vol. XI, núm. 27, abril, 2005, pp. 213-214, Universidad Tecnológica de Pereira-Colombia.

[2].Algunos tópicos en teoría de números: números Mersenne, teorema Dirichlet Números Fermat. Campo Elías González Pineda, Sandra Milena García. Scientia Et Technica, vol. XVI, núm. 48, agosto, 2011, pp. 185-190,Universidad Tecnológica de Pereira-Colombia.

[3] Una bella relación entre la conjetura de Goldbach y el teorema de Dirichlet. Scientia et technica 20 (3), 294-300. 2015/9/30Universidad Tecnológica de Pereira.

[4] https://arxiv.org/ftp/arxiv/papers/1610/1610.00598.pdf La maravillosa función y ecuación cuadrática. Campo Elías González Pineda.